\documentclass[a4paper, 10pt]{article}

\usepackage[T1]{fontenc}
\usepackage[latin1]{inputenc}
\usepackage{lmodern}  
\usepackage[english]{babel}
\usepackage{csquotes} 

\usepackage[a4paper]{geometry}
\usepackage{textcomp} 
\usepackage{graphicx} 
\usepackage{amsmath} 
\usepackage{amsthm} 
\usepackage{amssymb} 
\usepackage{amsfonts}
\usepackage{mathtools}
\usepackage{thmtools} 
\usepackage{mathrsfs} 
\usepackage{mathdots} 
\usepackage{wasysym} 
\usepackage{pdflscape} 
\usepackage{extpfeil} 
\usepackage{verbatim} 
\usepackage{longtable} 
\usepackage{bigstrut} 
\usepackage[table, dvipsnames]{xcolor} 
\usepackage{floatpag} 
\usepackage{enumitem} 

\usepackage{xparse} 
\usepackage{mathscinet} 
\usepackage{adjustbox} 
\usepackage[obeyDraft]{todonotes} 
\usepackage{pdfpages} 
\usepackage{caption}
\usepackage{subcaption} 
\usepackage{rotating}

\usepackage[final,                  
            colorlinks,             
            plainpages = false,     
            linktoc=all,            
            linkcolor=black,        
            citecolor=black,         
            urlcolor=black,          
            filecolor=black         
            ]{hyperref}

\usepackage[citestyle=alphabetic,   
            bibstyle=alphabetic,    
            maxbibnames=50,         
            maxcitenames=3,         
            autocite=inline,        
            block=space,            
            hyperref=true,          
            backref=false,           
            backrefstyle=three,     
            date=short,             
            arxiv=pdf,              
            isbn=false,             
            url=false,              
            doi=false,               
            eprint=true,            
            giveninits=true,        
            sorting=nyt,            
            backend=bibtex8,
            ]{biblatex}

\bibliography{billiards}

\hypersetup{pdftitle = {Correction of the Lusztig-Williamson Billiards Conjecture},
            pdfauthor = {Lars Thorge Jensen},
            pdfkeywords = {p-canonical basis} {Hecke algebra} {tilting characters} {SL3},
            pdfdisplaydoctitle      
            }
            
\usepackage[arrow, matrix, curve]{xy} 
\usepackage{tikz} 
\usetikzlibrary{arrows, arrows.meta, bending, fit, intersections, calc, external, shapes.misc, shapes.geometric, decorations.markings, decorations.pathreplacing}
\usepackage{pgfplots}
\pgfplotsset{compat=1.14}

\allowdisplaybreaks 

\setkeys{Gin}{draft=false}
            
\declaretheoremstyle[
spaceabove=\topsep, spacebelow=\topsep,
headfont=\normalfont\bfseries,
notefont=\mdseries, notebraces={(}{)},
bodyfont=\itshape,
postheadspace=\newline
]{break}

\declaretheoremstyle[
spaceabove=\topsep, spacebelow=\topsep,
headfont=\normalfont\bfseries,
notefont=\mdseries, notebraces={}{},
bodyfont=\itshape,
postheadspace=\newline
]{refbreak}

\declaretheorem[title=Theorem, style=plain, numberwithin=section]{thm}

\declaretheorem[title=Conjecture, style=break, numberlike=thm]{conj}

\declaretheorem[title=Remark, style=remark, numberlike=thm]{remark}

\declaretheorem[name=Heuristic, numbered=no]{heuristic}

\declaretheorem[name=Merge Rule, numbered=no]{mergeRule}

\usepackage[capitalize]{cleveref} 
                          
\crefname{thm}{Theorem}{Theorems}
\crefname{prop}{Proposition}{Propositions}
\crefname{lem}{Lemma}{Lemmata}
\crefname{cor}{Corollary}{Corollaries}
\crefname{rem}{Reminder}{Reminders}
\crefname{defn}{Definition}{Definitions}

\crefname{thmlab}{Theorem}{Theorems}
\crefname{proplab}{Proposition}{Propositions}
\crefname{lemlab}{Lemma}{Lemmata}
\crefname{corlab}{Corollary}{Corollaries}
\crefname{remlab}{Reminder}{Reminders}
\crefname{conj}{Conjecture}{Conjectures}

\crefname{thmreflab}{Theorem}{Theorems}
\crefname{propreflab}{Proposition}{Propositions}
\crefname{lemreflab}{Lemma}{Lemmata}
\crefname{correflab}{Corollary}{Corollaries}
\crefname{remreflab}{Reminder}{Reminders}
\crefname{conjref}{Conjecture}{Conjectures}

\crefname{remark}{Remark}{Remarks}
\crefname{claim}{Claim}{Claims}
\crefname{ex}{Example}{Examples}

\crefname{section}{Section}{Sections}
\crefname{figure}{Figure}{Figures}
\crefname{equation}{}{}
\crefname{ass}{Assumption}{Assumptions}

\CompileMatrices

\entrymodifiers={+!!<0pt,\fontdimen22\textfont2>} 

\def\clap#1{\hbox to 0pt{\hss#1\hss}}

\makeatletter
\def\underbracket{%
    \@ifnextchar[{\@underbracket}{\@underbracket [\@bracketheight]}%
}
\def\@underbracket[#1]{%
    \@ifnextchar[{\@under@bracket[#1]}{\@under@bracket[#1][0.4em]}%
}
\def\@under@bracket[#1][#2]#3{
    \mathop{\vtop{\m@th \ialign {##\crcr $\hfil \displaystyle {#3}\hfil $%
    \crcr \noalign {\kern 3\p@ \nointerlineskip }\upbracketfill {#1}{#2}
    \crcr \noalign {\kern 3\p@ }}}}\limits}

\def\upbracketfill#1#2{$\m@th \setbox \z@ \hbox {$\braceld$}
    \edef\@bracketheight{\the\ht\z@}\bracketend{#1}{#2}
    \leaders \vrule \@height #1 \@depth \z@ \hfill
    \leaders \vrule \@height #1 \@depth \z@ \hfill \bracketend{#1}{#2}$}

\def\bracketend#1#2{\vrule height #2 width #1\relax}
\makeatother

\makeatletter
\def\thmt@refnamewithcomma #1#2#3,#4,#5\@nil{%
  \@xa\def\csname\thmt@envname #1utorefname\endcsname{#3}%
  \ifcsname #2refname\endcsname
    \csname #2refname\expandafter\endcsname\expandafter{\thmt@envname}{#3}{#4}%
  \fi
}
\makeatother

\makeatletter
\newcommand*\rel@kern[1]{\kern#1\dimexpr\macc@kerna}
\newcommand*\widebar[1]{%
  \begingroup
  \def\mathaccent##1##2{%
    \rel@kern{0.8}%
    \overline{\rel@kern{-0.8}\macc@nucleus\rel@kern{0.2}}%
    \rel@kern{-0.2}%
  }%
  \macc@depth\@ne
  \let\math@bgroup\@empty \let\math@egroup\macc@set@skewchar
  \mathsurround\z@ \frozen@everymath{\mathgroup\macc@group\relax}%
  \macc@set@skewchar\relax
  \let\mathaccentV\macc@nested@a
  \macc@nested@a\relax111{#1}%
  \endgroup
}
\makeatother

\makeatletter
\newcommand{\subjclass}[2][1991]{%
  \let\@oldtitle\@title%
  \gdef\@title{\@oldtitle\footnotetext{#1 \emph{Mathematics subject classification.} #2}}%
}
\newcommand{\keywords}[1]{%
  \let\@@oldtitle\@title%
  \gdef\@title{\@@oldtitle\footnotetext{\emph{Key words and phrases.} #1.}}%
}
\makeatother

\makeatletter
\newcommand{\extp}{\@ifnextchar^\@extp{\@extp^{\,}}}
   \def\@extp^#1{\mathop{\bigwedge\nolimits^{\!#1}}}
\makeatother

\def\un{\underline}

\DeclareMathOperator{\Hom}{Hom}
\DeclareMathOperator{\End}{End}

\newcommand{\id}{\mathrm{id}}
\DeclareMathOperator{\sgn}{sgn} 

\newcommand{\defeq}{\ensuremath{\coloneqq}}

\newcommand{\Z}{\ensuremath{\mathbb{Z}}}
\newcommand{\N}{\ensuremath{\mathbb{N}}}

\newcommand{\R}{\ensuremath{\mathbb{R}}}

\newcommand{\heck}{\ensuremath{H}}

\newcommand{\asphMod}[1][]{\ensuremath{AS^{#1}_v}}

\newcommand{\std}[2][n]{\ensuremath{#1_{#2}}}
\newcommand{\kl}[2][n]{\ensuremath{\underline{#1}_{#2}}}
\DeclareDocumentCommand{\pkl}{O{p} O{n} m}{\ensuremath{\prescript{#1}{}{\underline{#2}}_{#3}}}
\DeclareDocumentCommand{\pklgen}{O{p} O{n} m m}{\ensuremath{\prescript{#1}{}{\underline{#2}}^{#4}_{#3}}}

\newcommand{\desc}[1]{\ensuremath{\mathcal{#1}}}

\newcommand{\rt}[1][]{\ensuremath{\alpha_{#1}}}
\newcommand{\cort}[1][]{\ensuremath{\alpha_{#1}^{\vee}}}
\newcommand{\rts}{\ensuremath{\Phi}}

\newcommand{\rtbasis}{\ensuremath{\Sigma}}

\newcommand{\charlat}{\ensuremath{\mathscr{X}}}

\newcommand{\fin}{\ensuremath{\text{f}}}

\newcommand{\finWeylGrp}{\ensuremath{W_{\fin}}}
\newcommand{\finRefl}{\ensuremath{S_{\fin}}}
\newcommand{\affWeylGrp}[1][]{\ensuremath{W}_{#1}}
\newcommand{\affRefl}{\ensuremath{S}}

\newcommand{\fW}{\ensuremath{{}^{\fin}\affWeylGrp}}

\newcommand{\dotAct}[1][p]{\ensuremath{\cdot_{#1}}}

\newcommand{\Homnl}[2][]{
   \ifthenelse{ \equal{#1}{} } { \ensuremath{\Hom_{\nless #2}} } 
   { \ensuremath{\Hom_{\nless #2, #1}} } }
\newcommand{\Endnl}[2][]{
   \ifthenelse{ \equal{#1}{} } { \ensuremath{\End_{\nless #2}} } 
   { \ensuremath{\End_{\nless #2, #1}} } }

\newcommand{\labels}{\ensuremath{\mathcal{M}}}
\newcommand{\gWall}{\ensuremath{\Gamma_{\text{wall}}}}

\newcommand{\Address}{
  \bigskip{\footnotesize
  \textsc{\'Ecole polytechnique f\'ed\'erale de Lausanne, Route Cantonale, CH-1015 Lausanne, Switzerland}\par\nopagebreak
   \textsc{Institute for Advanced Study, 1 Einstein Drive, Princeton, NJ 08540, USA}\par\nopagebreak
  \textit{E-mail address}, Lars~Thorge~Jensen: \texttt{ltjensen@ias.edu}
}}

\tikzset{%
  DynNode/.style={circle, inner sep=2pt, draw=black, fill=white},
  Greater/.style={pos=0.65, inner sep=0mm, outer sep=0mm},
  highlight/.style={rectangle,rounded corners,fill=red!15,draw=red,
    fill opacity=0.5,thick},
  root/.style={draw, color=black, thick, ->},
  plane/.style={draw, color=black, very thin},
  origin/.style={fill, color=black},
  sline/.style={color=Red, thick},
  tline/.style={color=NavyBlue, thick},
  uline/.style={color=Goldenrod, thick},
  bendBelow/.style={bend left=70, looseness=2},
  bendAbove/.style={bend right=70, looseness=2},
  object/.style={circle, fill, inner sep=1.5pt, outer sep=0mm},
  labelling/.style={outer sep=0mm, inner sep=0mm},
  1morph/.style={->, shorten >= 0.5pt, >=stealth'},
  2morph/.style={-implies,double,double equal sign distance,
                 shorten >=2pt, shorten <=3pt},
  spot/.style={color=black, thin, dashed},
  s/.style={color=Red},
  t/.style={color=NavyBlue},
  u/.style={color=Goldenrod},
  line/.style={draw, line width=2pt},
  dot/.style={fill, thin},
  sph/.style={fill, color=black!20, opacity=0.5},
  squig/.style={decoration=snake, decorate, ->},
  root/.style={draw, color=black, line width=2pt, ->},
  myptr/.style={decoration={markings,mark=at position 1 with %
    {\arrow[scale=3,>=stealth]{>}}},postaction={decorate}},
  on each segment/.style={
    decorate,
    decoration={
      show path construction,
      moveto code={},
      lineto code={
        \path [#1]
        (\tikzinputsegmentfirst) -- (\tikzinputsegmentlast);
      },
      curveto code={
        \path [#1] (\tikzinputsegmentfirst)
        .. controls
        (\tikzinputsegmentsupporta) and (\tikzinputsegmentsupportb)
        ..
        (\tikzinputsegmentlast);
      },
      closepath code={
        \path [#1]
        (\tikzinputsegmentfirst) -- (\tikzinputsegmentlast);
      },
    },
  },
  mid arrow/.style={postaction={decorate,decoration={
        markings,
        mark=at position .5 with {\arrow[#1]{stealth}}
      }}},
  sline/.style={draw, line width=1pt, postaction={on each segment={mid arrow=black}}},
  sregion/.style={fill, opacity=0.2},
  st/.style={fill=Fuchsia},
  su/.style={fill=YellowOrange},
  tu/.style={fill=ForestGreen},
  clabel/.style={fill=none, red}, 
  str/.style={<->}
}

\begin{document}

\title{Correction of the Lusztig-Williamson Billiards Conjecture}
\author{Lars Thorge Jensen}
\subjclass[2000]{Primary 20C20, 17B10, 20C30}
\keywords{tilting modules; billiards; p-canonical basis; p-Kazhdan-Lusztig basis; symmetric group}
\date{} 

\maketitle

\begin{abstract}
   A new algorithm allows us to calculate many new tilting characters for $SL_3$, $SP_4$,
   $G_2$, $SL_4$ and potentially many other groups. These calculations show that the 
   Lusztig-Williamson Billiards Conjecture needs to be corrected. In this paper we 
   present the new results calculated for $SL_3$ and a correction of the conjecture.
\end{abstract}

\tableofcontents

\section{Introduction}

Recent results show that the characters of indecomposable tilting modules of reductive
algebraic groups in characteristic $p$ are given by the $p$-Kazhdan-Lusztig basis of 
the corresponding anti-spherical module (see \cite{RWTiltPCan, AMRWTiltChars,
RWSmithTreumann, BRAction}) for all primes $p$. In joint work with Geordie Williamson, 
the author has developed and implemented (in Magma, see \cite{Magma}) a new algorithm 
to calculate the $p$-Kazhdan-Lusztig basis of the anti-spherical module. This 
algorithm will be described in detail in a separate paper. Even though the 
calculations which the paper \cite{LuWTiltCharsSL3} is based on can now be carried 
out in less than 24 hours, this should not be considered a full solution to the problem.

The Lusztig-Williamson Billiards Conjecture gives a combinatorial description
of second generation tilting characters for $SL_3$ in terms of billiard balls bouncing
around in the $p$-alcove geometry of the dominant cone. For background
on the generational philosophy we refer the reader to \cite[\S2]{LuWTiltCharsSL3} and
the references therein. The new calculations show that the Lusztig-Williamson 
Billiards Conjecture actually predicts second generation tilting characters that are 
``too big'' and thus needs to be corrected.\footnote{The data Lusztig and Williamson could calculate at the time did not show the phenomena.} We propose a correction to the conjecture that
is coherent with all the data we can currently calculate for $SL_3$. 

\subsection{Structure of the Paper}

\begin{description}
   
   \item[\Cref{secNot}] We will fix important notation and recall generations for
      the $p$-Kazhdan-Lusztig basis of the anti-spherical module.
   \item[\Cref{secConj}] In this section, we will explain the correction to the
      Lusztig-Williamson Billiards Conjecture. The reader will understand why
      the original conjecture predicts anti-spherical $p$-Kazhdan-Lusztig basis
      elements that are ``too big''.
   \item[\Cref{secData}] Finally, we will explain how to obtain the figures from
      the anti-spherical $p$-Kazhdan-Lusztig basis calculated by our algorithm
      and mention all the examples for the newly proposed merging rule that can
      be found in the current data.
\end{description}         
         
\subsection{Acknowledgements}
The counterexample to the Lusztig-Williamson Billiards Conjecture was discovered while
the author was trying to understand the combinatorics of the Smith-Treumann localization
of the anti-spherical module during a research at the Sydney Mathematical Research Institute.
For that reason, the author would like to thank Geordie Williamson for the invitation
and the SMRI for the excellent working conditions. In addition, the author would like
to thank Allan Steel for helping to improve the runtime of the algorithm.
During that period the author has received funding from the European Research Council 
(ERC) under the European Union's Horizon 2020 research and innovation programme 
(grant agreement No 677147).

\section{Notation}
\label{secNot}

Let $G$ be a split simple and simply connected algebraic group over an algebraically
closed field of characteristic $p$ and fix a Borel subgroup together with a maximal
torus $T \subseteq B \subseteq G$. We will imitate the notation of 
\cite[\S 2]{LuWTiltCharsSL3} to make the transition between the two papers as easy
as possible for the reader:
\begin{align*}
   \charlat, \charlat_+, \charlat_{++}\text{:} \quad&\text{weights, dominant weights, strictly dominant weights;}\\
   \rts, \rts_+\text{:} \quad&\text{roots, positive roots;}\\
   \rtbasis\text{:} \quad&\text{simple roots;}\\
   \rho, \cort[0]\text{:} \quad&\text{half-sum of $\rts_+$, highest short coroot;}\\
   \finWeylGrp, \affWeylGrp\text{:} \quad&\text{finite Weyl group, affine Weyl group;}\\
   \finRefl, \affRefl\text{:} \quad&\text{finite simple reflections, affine simple reflections;}\\
   \fW\text{:} \quad&\text{minimal coset representatives for the cosets $\finWeylGrp 
               \backslash \affWeylGrp$;}\\
   \dotAct, h\text{:} \quad&\text{$p$-dilated dot action of $\affWeylGrp$ on $\charlat$, Coxeter number.}
\end{align*}

Let $\heck$ (resp. $\heck_{\fin}$) be the Iwahori-Hecke algebra of the affine (resp. finite)
Weyl group $\affWeylGrp$ (resp. $\finWeylGrp$) over $\Z[v, v^{-1}]$ and consider the
anti-spherical module:
\[ \asphMod \defeq \sgn_v \otimes_{\heck_{\fin}} \heck = \bigoplus_{x\in \fW} \Z[v, v^{-1}] \std{x} \]
The anti-spherical module has a Kazhdan-Lusztig basis $\{\kl{x} \; \vert \; x \in \fW\}$
and a $p$-Kazhdan-Lusztig basis $\{\pkl{x} \; \vert \; x \in \fW\}$. Lusztig and 
Williamson expect that for all $p > h$ and $x \in \fW$ there exist
elements $\pklgen{x}{1}$, $\pklgen{x}{2}$, $\dots$, $\pklgen{x}{\infty} \in \asphMod$ 
(see \cite[\S 2.4]{LuWTiltCharsSL3}) such that:
\begin{enumerate}
   \item $\pklgen{x}{1} = \kl{x}$;
   \item $\pklgen{x}{\infty} = \pkl{x}$;
   \item $\pklgen{x}{n}$ is a $\Z_{\geqslant 0}[v, v^{-1}]$-linear combination of
         $\{ \pklgen{y}{n-1} \; \vert \; y \in \fW\}$ for all $n \geqslant 1$;
   \item if $\langle \cort[0], x \dotAct 0 + \rho \rangle \leqslant p^{n+1}$ then
         $\pklgen{x}{n} = \pklgen{x}{n+1} = \dots = \pklgen{x}{\infty}$.
\end{enumerate}
This expectation was motivated by the generational philosophy and based on the 
hypothesis that the tilting character generations can be lifted to the non-trivial 
grading of the category of tilting modules given by the anti-spherical category 
(see \cite[Theorem 5.3]{RWTiltPCan}).

\section{The Corrected Conjecture}
\label{secConj}

The Lusztig-Williamson Billiards Conjecture gives a combinatorial description
of $\{ \pklgen{x}{2} \; \vert \; x \in \fW\}$ for $G = SL_3$. For that reason,
we will fix $G = SL_3$ for the rest of the paper. Thus, we have 
$\charlat = \Z \varpi_1 \oplus \Z \varpi_2$, 
$\charlat_{+} \defeq \Z_{\geqslant 0} \varpi_1 \oplus \Z_{\geqslant 0} \varpi_2$
$\charlat_{++} \defeq \Z_{>0} \varpi_1 \oplus \Z_{>0} \varpi_2$, etc..

Fix $\ell \leqslant 3$. Note that $\ell$ will be prime for the representation
theoretic applications. Denote by 
$\labels \defeq \Z_{\geqslant 0} \times \{v^k \; \vert \; k \in \Z\}$
the set of \emph{labels}. The element $(m, v^k) \in \labels$ will be
written as $m(v^k)$. A \emph{labelled point} is an element of
$\charlat \times \labels$. Throughout this section all set-theoretic operations 
(i.e. unions, differences, \dots) should be understood in the context of 
multisets (i.e. sets with multiplicities).

In \cite[\S4]{LuWTiltCharsSL3} Lusztig-Williamson describe an algorithm 
which in three steps produces a multiset $\widetilde{Z}$ from which one can conjecturally obtain the
set $\{ \pklgen{x_i}{2} \; \vert \; i \geqslant 0\}$ where
\[ x_0 \defeq \id\text{, }x_1 \defeq s_0\text{, }x_2 \defeq s_0 s_1\text{, }
x_3 \defeq s_0 s_1 s_2\text{, }x_4 \defeq s_0 s_1 s_2  s_0, \dots \]
and $\{s_0\} = \affRefl \setminus \finRefl$ is the affine reflection.

In the correction, the first and the third step of the algorithm remain unchanged
and we refer the reader to \cite[\S4.2 and \S4.4]{LuWTiltCharsSL3} for their
description. The first step starts with the labelled point $(0, 0(v^0))$ and
produces a set $X$ by extending along a wall of the dominant cone in the direction
of $\varpi_1$. In each step, certain labelled points are designated as seeds
and serve as input for the next step. The seeds in $X$ are the labelled points
of the form $(k\ell\varpi_1, 2k\ell(v^0))$ for $k \geqslant 1$.

Before we explain the correction of the second step, we should mention
how to put everything together in the end. The third step of the algorithm
takes the multiset $Y$ produced by the second step and extends it within
the interior of each $\ell$-alcove to produce a multiset $Z$. Finally, 
we set $\widetilde{Z} \defeq Z \setminus X$.

In order to formulate the conjecture, we will need some more notation
from \cite[\S6]{LuWTiltCharsSL3}. Recall our prime $p$ from \cref{secNot}
and consider the multiset $\widetilde{Z}$ defined as above for $\ell = p$.
Consider the $\Z$-linear map $\varphi: \Z[v] \rightarrow \Z[v, v^{-1}]$
defined via $v^0 \mapsto 1$ and $v^k \mapsto v^k + v^{-k}$ for $k > 0$.
For an element $w \in \affWeylGrp$ we will denote by $\desc{R}(w) \defeq 
\{ s \in \affRefl \; \vert \; ws < w \}$ its right descent set.
Fix the fundamental alcove
\[ \mathcal{A}_0 \defeq \{ \lambda \in \charlat \otimes_{\Z} \R \; \vert \; 
   0 < \langle \cort, \lambda \rangle < 1 \text{ for all } \rt \in \rts_+ \} \]
and recall the action of $\affWeylGrp$ on the set of alcoves induced by
the continuous action of $\affWeylGrp$ on $\charlat \otimes_{\Z} \R$.
As described in \cite[\S6]{LuWTiltCharsSL3}, for $\mu \in \charlat_{++}$ and 
$s\in \affRefl$ there exists a unique element $x_{\mu}^s \in \fW$ such that 
$s \in \desc{R}(x_{\mu}^s)$ and the open box 
\[ B_{\mu} \defeq \{ \lambda \in \charlat \otimes_{\Z} \R \; \vert \; 
   \langle \cort, \mu \rangle < \langle \cort, \lambda \rangle < 
   \langle \cort, \mu\rangle + 1 \text{ for all } \rt \in \rtbasis \} \]
contains the alcove $x\mathcal{A}_0$.

Define ${}^p \zeta_0 \defeq \kl{x_0}$. For $i > 0$ let $s$ be the unique
simple reflection in $\desc{R}(x_i)$ and define the element
\[ {}^p \zeta_i \defeq \kl{x_i} + \sum_{\substack{(\mu, n(v^k)) \in \widetilde{Z}\\
   n \in \{i, i-1, i-2\}}} \varphi(v^k) \kl{x_{\mu}^s}\text{.}\]
\begin{conj}[Corrected Lusztig-Williamson Billiards Conjecture]
   We have:
   \begin{enumerate}
      \item ${}^p \zeta_i = \pkl{x_i}$ for $0 \leqslant i < 2p(p+1)$,
      \item ${}^p \zeta_i = \pklgen{x_i}{2}$ \text{ for all } $0 \leqslant i$.
   \end{enumerate}
\end{conj}

\subsection{Corrected 2. Step: Dynamics on the Walls}

In order to describe the second step of the algorithm, we first need to recall
some definitions from \cite[\S4.3]{LuWTiltCharsSL3} and we will refer the
reader to the original source for some beautiful illustrations of these
definitions.

\emph{Corner points} are points $\lambda \in \charlat_+$ such that 
$\langle \lambda, \cort \rangle \in \ell\Z$ for all $\rt \in \rts_+$. Consider the 
directed graph $\Gamma$ on the vertex set $\charlat_+$ with edges 
$\lambda \rightarrow \lambda + \gamma$ for $\lambda, \lambda + \gamma \in \charlat_+$ 
and $\gamma \in \{\varpi_1, \varpi_2 - \varpi_1, -\varpi_2\}$.
A point $\mu \in \charlat_+$ is called an \emph{almost corner} if there
exists a corner point $\lambda \in \charlat_+$ and an edge $\lambda \rightarrow \mu$
in $\Gamma$.

For the dynamics on the walls, the following subgraphs $\gWall$ and $\Gamma_{\ell}$ 
of $\Gamma$ will be important. The vertex set of $\Gamma_{\ell}$ consists of weights 
$\lambda \in \charlat_+$ such that $\langle \lambda, \cort \rangle \in \ell\Z$ for 
some $\rt \in \rts_+$. The only edges in $\Gamma_{\ell}$ are those $\lambda \rightarrow
\lambda'$ where $\langle \lambda, \cort \rangle = \langle \lambda', 
\cort \rangle \in \ell\Z$ for some $\rt \in \charlat_+$. Define $\gWall$ 
to be the induced subgraph of $\Gamma_{\ell}$ where we remove all the 
vertices that lie along the walls of the dominant
cone except for $\{ k\ell\varpi_1, k\ell\varpi_2 \; \vert \; k \geqslant 1\}$.
Recall that the labelled points of the form $(k\ell\varpi_1, 2k\ell(v^0))$ are designated
as \emph{seeds} in the first step of the algorithm.

Fix a labelled point $(\mu, n(v^k))$ with $\mu \in \gWall$ such that there is a 
unique edge with source $\mu$ in $\gWall$. First, Lusztig and Williamson define 
three elementary operations to obtain labelled points in $\gWall \times \labels$:
\begin{enumerate}
   \item A \emph{rest} gives the labelled point $(\mu, (n+3)(v^{k+1}))$.
   \item A \emph{small step} gives the labelled point $(\mu', (n+2)(v^k))$ where
         $\mu \rightarrow \mu'$ is the unique edge in $\gWall$ starting in $\mu$.
   \item If $\mu$ is neither a corner nor an almost corner point, a \emph{giant leap}
         gives either one or two new labelled points as follows: Let $d$ be the
         direction of the unique arrow in $\gWall$ starting in $\mu$. First take
         $j < \ell-1$ steps in the graph $\gWall$ in direction $d$ until we reach
         a corner point, say $\lambda$. Then continue with another $\ell - 1 - j$ 
         steps in all directions different from $d$ from $\lambda$. (There are either
         one or two such directions depending on whether we are close to
         the walls of the dominant cone or not.) The giant leap consists of
         the resulting points labelled by $(n + 2\ell + 1)(v^{k+1})$.
         In the following we will also call this a 
         \emph{giant leap towards $\lambda$}.
\end{enumerate}

Then they define one iteration of the algorithm producing a multiset of 
elements in $\gWall \times \labels$. Some of these labelled points
will be designed as seeds and thus serve as input for the third step of the
algorithm. Using the fixed labelled point $(\mu, n(v^k))$ as input,
we proceed as follows:
\begin{enumerate}
   \item If $\mu$ is a \emph{corner point}, then $\ell$ new labelled points are obtained by
      taking $\ell-1$ small steps to obtain $m_1, m_2, \dots, m_{\ell-1}$ and then a rest
      to produce $m_{\ell}$. The last labelled point $m_{\ell}$ is designated as a seed.
      This operation is called \emph{resting once}.
   \item If $\mu$ is an \emph{almost corner point}, then $\ell$ new labelled points are
      obtained by taking first a rest to produce $m_1$, then $\ell-2$ small steps to obtain
      $m_2, \dots m_{\ell-1}$ and finally another rest to produce $m_{\ell}$. The last
      labelled point $m_{\ell}$ is designated as a seed. We will call this operation 
      \emph{resting twice}.
   \item  If $\mu$ is neither a \emph{corner} nor an \emph{almost corner point}, then
      we take a giant leap to produce one or to new labelled points, each of which
      are designated as seeds.
\end{enumerate}

Given a seed $q_k = (k\ell\varpi_1, 2k\ell(v^0))$ in $X$ for some $k \geqslant 1$, 
we iterate this process as follows:
Starting with $Q_0 \defeq \{q_k\}$ we obtain a sequence of multisets $Q_1, Q_2, \dots$ 
where $Q_i$ is the multiset obtained by taking the union of the outputs of the 
above procedure applied to each seed in $Q_{i-1}$  
and then applying the following geometric merge rule for every corner point $\lambda$:
\begin{mergeRule}
   If giant leaps towards $\lambda$ lead to a superposition of labelled points
   (i.e. at least one labelled point occuring with multiplicity two), then among 
   the seeds produced by these giant leaps towards $\lambda$ we only keep the 
   superposed seeds and reduce their multiplicity to one.
\end{mergeRule}
Finally we say that the union $Y_k = \bigcup_{i \geqslant 1} Q_i$ is the result
of \emph{applying dynamics on the walls to $q_k$}. Note that the geometric 
merge rule is the only new ingredient in the correction.

The output of the corrected second step is given by the multiset
\[ Y \defeq X \cup \bigcup_{k \geqslant 1} Y_k \]
together with the information which elements of $Y$ are designated as seeds.

In order to make the geometric merge rule more explicit, we want to explain the 
three possible cases for every corner point $\lambda$:
\begin{enumerate}
   \item If there is one giant leap towards $\lambda$, it produces one or two 
      labelled points (depending on whether $\lambda$ is close to the walls
      of the dominant cone or not).
   \item If there are two giant leaps towards $\lambda$, they produce one labelled 
     point (where the coefficient is half the sum of the superposed coefficients).
     We will call this case \emph{merging of type \RN{2} around $\lambda$}.
   \item If there are three giant leaps towards $\lambda$, they produce three labelled
      points (each with half the sum of the superposed coefficients as coefficient). 
      This case is called \emph{merging of type \RN{3} around $\lambda$}.
\end{enumerate}

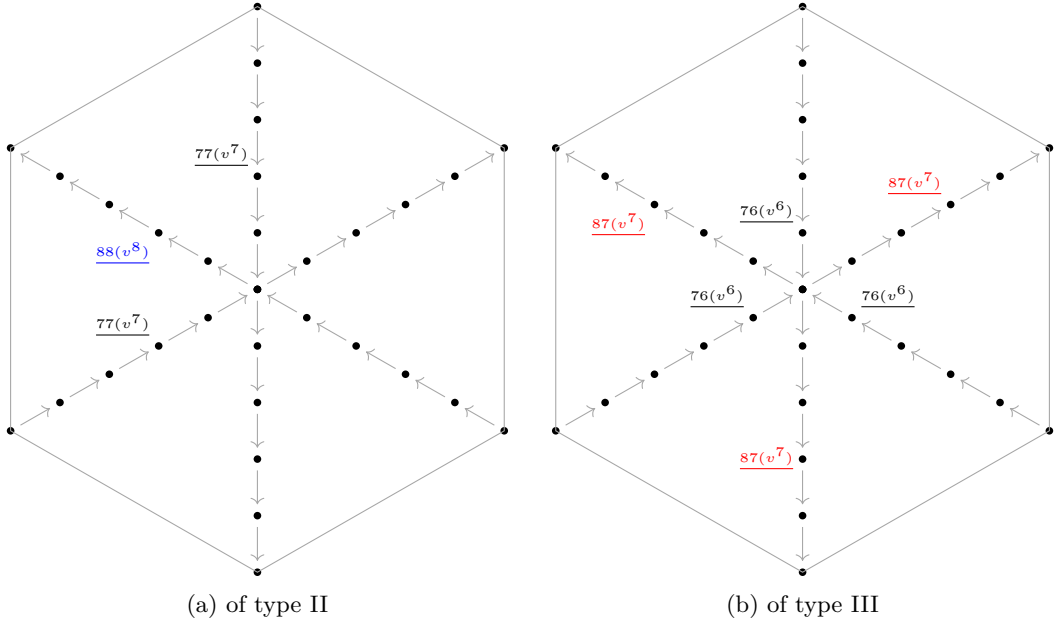
\begin{figure}[htb]
   \centering
   
   \begin{subfigure}{.48\textwidth}
      \centering
      \begin{tikzpicture}[scale=0.75, rotate=30, point/.style={inner sep=0pt,outer sep = 1mm,circle,fill,minimum size=1mm}]
         \pgftransformcm{1}{0}{0.5}{sqrt(3)/2}{\pgfpoint{0cm}{0cm}};   
         \foreach \y in {0,1, ..., 10} 
         {
            \node (x5-y\y) at (5,\y) [point] {};
            \ifthenelse{\y > 0}{
               \pgfmathtruncatemacro{\yprev}{\y - 1}
               \draw[color=gray!70,->] (x5-y\y) -- (x5-y\yprev);}
               {}
         }
         \foreach \y in {0, ..., 10}
         {
            \pgfmathtruncatemacro{\x}{10 - \y}
            \node (x\x-y\y) at (\x,\y) [point] {};
            \ifthenelse{\y > 0}{
               \pgfmathtruncatemacro{\xprev}{\x+1}
               \pgfmathtruncatemacro{\yprev}{\y-1}
               \draw[color=gray!70,->] (x\xprev-y\yprev) -- (x\x-y\y);
            }{}
         }
         \foreach \x in {0, ..., 10} 
         {
            \node (x\x-y5) at (\x,5) [point] {};
            \ifthenelse{\x > 0}{
               \pgfmathtruncatemacro{\xprev}{\x-1}
               \draw[color=gray!70,->] (x\xprev-y5) -- (x\x-y5);
            }{}
         }
         \draw[color=gray!70,->] (5,0) -- (10,0) -- (10,5) -- (5,10) -- (0,10) -- (0,5) -- cycle;
         \node[above left] at (5,7) {\tiny $\un{77(v^7)}$};  
         \node[above left] at (3,5) {\tiny $\un{77(v^7)}$};  
         
         \node[below left] at (3,7) {\color{blue}\tiny $\un{88(v^8)}$};  
      \end{tikzpicture}
      
      \caption{of type \RN{2}}
      \label{figMergingOfTypeII}
   \end{subfigure}
   \begin{subfigure}{.48\textwidth}
      \centering
      \begin{tikzpicture}[scale=0.75, rotate=30, point/.style={inner sep=0pt,outer sep = 1mm,circle,fill,minimum size=1mm}]
         \pgftransformcm{1}{0}{0.5}{sqrt(3)/2}{\pgfpoint{0cm}{0cm}};   
         \foreach \y in {0,1, ..., 10} 
         {
            \node (x5-y\y) at (5,\y) [point] {};
            \ifthenelse{\y > 0}{
               \pgfmathtruncatemacro{\yprev}{\y - 1}
               \draw[color=gray!70,->] (x5-y\y) -- (x5-y\yprev);}
               {}
         }
         \foreach \y in {0, ..., 10}
         {
            \pgfmathtruncatemacro{\x}{10 - \y}
            \node (x\x-y\y) at (\x,\y) [point] {};
            \ifthenelse{\y > 0}{
               \pgfmathtruncatemacro{\xprev}{\x+1}
               \pgfmathtruncatemacro{\yprev}{\y-1}
               \draw[color=gray!70,->] (x\xprev-y\yprev) -- (x\x-y\y);
            }{}
         }
         \foreach \x in {0, ..., 10} 
         {
            \node (x\x-y5) at (\x,5) [point] {};
            \ifthenelse{\x > 0}{
               \pgfmathtruncatemacro{\xprev}{\x-1}
               \draw[color=gray!70,->] (x\xprev-y5) -- (x\x-y5);
            }{}
         }
         \draw[color=gray!70,->] (5,0) -- (10,0) -- (10,5) -- (5,10) -- (0,10) -- (0,5) -- cycle;
         \node[above left] at (5,6) {\tiny $\un{76(v^6)}$};  
         \node[above right] at (6,4) {\tiny $\un{76(v^6)}$};  
         \node[above left] at (4,5) {\tiny $\un{76(v^6)}$};  
         
         \node[below left] at (2,8) {\color{red}\tiny $\un{87(v^7)}$};  
         \node[above left] at (8,5) {\color{red}\tiny $\un{87(v^7)}$};  
         \node[left] at (5,2) {\color{red}\tiny $\un{87(v^7)}$};  
      \end{tikzpicture}
      
      \caption{of type \RN{3}}
      \label{figMergingOfTypeIII}
   \end{subfigure}
    \caption{Examples of Merging for  $\ell = 5$}
\end{figure}

\begin{remark}
   \begin{enumerate}
   \item Lusztig and Williamson were not aware of the merging phenomenon as
         they did not have a single occurrence of it in the data they could
         calculate.
   \item The merge rule annihilates \cite[Remark 4.3]{LuWTiltCharsSL3}. In other
         words, it prevents exponentially growing coefficients and thus the
         corrected conjecture no longer implies that decomposition numbers
         for symmetric groups display exponential growth.
   \end{enumerate}
\end{remark}

\subsection{A Larger Example and Some Observations}

Figures \ref{figl5k1} and \ref{figl5k2} show larger examples to 
illustrate the geometric merge rule. They show for $\ell = 5$ the multiset after 
the first ten iterations of the algorithm described in the second step 
applied to the seeds $(k\ell\varpi_1, 2k\ell(v^0))$ for $k \in \{1, 2\}$.
In this example, we have underlined all seeds and used the colours {\color{red}red}
(resp. {\color{blue}blue}) to mark seeds that arise from a merging of type
\RN{3} (resp. of type \RN{2}). When looking at this example the reader should 
observe the following:

\begin{itemize}
   \item All the seeds in $Q_j$ are obtained by applying the same type of operation to
         some seed in $Q_{j-1}$. Looking at this operation for the sequence of
         multisets $Q_1, Q_2, Q_3, \dots$, we get the following $\ell$-periodic pattern 
         of operations:
         \begin{align*}
               &\text{resting once, }\underbrace{
               \text{giant leap, \dots, giant leap}}_{(\ell-2)\text{-times}}\text{, resting twice,}\\
               &\overline{\underbrace{
               \text{giant leap, \dots, giant leap}}_{(\ell-1)\text{-times}}\text{, resting twice}}\text{, \dots}
         \end{align*}
         
   \item There are no seeds with a label $n(v^{k\ell})$ for $n\in \N$ and $k \geqslant 1$
         (simply because resting twice augments the exponent of $v$ by $2$ whereas
         all other operations increase it by $1$).
   \item For a labelled point $(\mu, n(v^k)) \in \gWall \times \labels$ we always have:
         \[ n \equiv k \mod \ell \]
   \item For a labelled point $(\mu, n(v^k)) \in \gWall \times \labels$,  the 
         residue of $n$ modulo $\ell$ determines the position of the labelled point 
         along the wall on which it lies.
   \item The seeds in $Q_j$ obtained from 
         any seed in $X$ are all of the form 
         \[ (\mu, n(v^{j + \lfloor \frac{j-1}{\ell -1} \rfloor}))\] 
         for various $\mu \in \charlat_+$. In other words, the exponent of $v$ 
         determines how many iterations are needed to produce this point.
   \item The previous observations together imply that there are no superposition
         of labelled points in $Y_j$ and $Y_k$ for $j \ne k$.
   \item All labelled points $(\mu, n(v^k))$ lying in $Y_k$ (i.e. the ones obtained 
         by applying the iterative algorithm in the second step to 
         $(k\ell \varpi_1, 2k\ell(v^0))$) satisfy $\mu = a \varpi_1 + b \varpi_2$ with
         \[ 0 < a < k\ell \quad \text{ and } \quad  0 < b \text{.}\]
         In other words, they lie in the first $k$ ``columns'' of the dominant cone
         (see Figures \ref{figl5k1} and \ref{figl5k2}).
   \item The keen reader might have noticed that the ``column'' depicted in Figure 
         \ref{figl5k1} matches up to the global label shift $n(v^k) \mapsto (n+10)(v^k)$
         the second ``column'' in Figure \ref{figl5k2} (without the left wall!).
         More generally, the mapping $Y_j \rightarrow Y_k$ for $j < k$ given by
         $(\mu, n(v^i)) \mapsto (\mu + (k-j)\ell\varpi_1, (n+2(k-j)\ell)(v^i))$
         is well-defined, injective and surjective onto the set
         \[  \{ (a \varpi_1, + b\varpi_2, n(v^i)) \in Y_k \; \vert \; (k-j) \ell < a \}
         \text{.} \]
         This suggests a striking self-similarity of the data!
         Moreover, it implies that for the second step it is enough to understand the
         combinatorics of the wall dynamics applied to the seed 
         $(k\ell\varpi_1, 2k\ell(v^0))$ in the limit for $k \rightarrow \infty$.
         From this, one can obtain the multiset $Y$ by taking unions of preimages
         of various maps of the type described above.
\end{itemize}
 
\begin{figure*}[p]
   \centering
   \begin{minipage}[b]{.35\textwidth}
      \centering
      \resizebox{\textwidth}{!}{
      \begin{tikzpicture}[scale=1.1, rotate=30, point/.style={inner sep=0pt,outer sep = 1mm,circle,fill,minimum size=1mm}]
         \pgftransformcm{1}{0}{0.5}{sqrt(3)/2}{\pgfpoint{0cm}{0cm}};   
         \draw[opacity=0.5, ->] (0,0) -- (5.5,0);
         \draw[opacity=0.5, ->] (0,0) -- (0,22);
         \foreach \y in {0,1, ..., 20} 
         {
            \node (x5-y\y) at (5,\y) [point] {};
            \ifthenelse{\y > 0}{
               \pgfmathtruncatemacro{\yprev}{\y - 1}
               \draw[color=gray!70,->] (x5-y\y) -- (x5-y\yprev);}
               {}
         }
         \foreach \t in {1, ..., 5}
         {
            \pgfmathtruncatemacro{\x}{5 - \t}
            \node (x\x-y\t) at (\x,\t) [point] {};         
            \pgfmathtruncatemacro{\xprev}{\x+1}
            \pgfmathtruncatemacro{\yprev}{\t-1}
            \draw[color=gray!70,->] (x\xprev-y\yprev) -- (x\x-y\t);
         }
         \foreach \ty in {0, 5, 10} {    
            \foreach \tx in {6, ..., 10}
            {
               \pgfmathtruncatemacro{\x}{10 - \tx}
               \pgfmathtruncatemacro{\y}{\ty + \tx}
               \node (x\x-y\y) at (\x,\y) [point] {};
               \pgfmathtruncatemacro{\xprev}{\x+1}
               \pgfmathtruncatemacro{\yprev}{\y-1}
               \draw[color=gray!70,->] (x\xprev-y\yprev) -- (x\x-y\y);
            }
         }
         \foreach \y in {5, 10, 15, 20} {
            \foreach \x in {1, ..., 5} 
            {
               \node (x\x-y\y) at (\x,\y) [point] {};
               \pgfmathtruncatemacro{\xprev}{\x-1}
               \draw[color=gray!70,->] (x\xprev-y\y) -- (x\x-y\y);
            }
         }
         
         \node[below right] at (5,0) {\tiny $\un{10(v^0)}$};
         \node[below left] at (4,1) {\tiny $12(v^0)$};
         \node[below left] at (3,2) {\tiny $14(v^0)$};
         \node[below left] at (2,3) {\tiny $16(v^0)$};
         \node[below left] at (1,4) {\tiny $18(v^0)$};
         \node[left] at (1,4) {\tiny $\un{21(v^1)}$};
         \node[above left] at (3,5) {\tiny $\un{32(v^2)}$};      
         \node[right] at (5,3) {\tiny $\un{43(v^3)}$};
         \node[below left] at (3,7) {\tiny $\un{43(v^3)}$};
         
         \node[above right] at (4,1) {\tiny $\un{54(v^4)}$};
         \node[label={[above right, yshift=1ex]:\tiny $57(v^5)$}] at (4,1) {};
         \node[above right] at (3,2) {\tiny $59(v^5)$};
         \node[above right] at (2,3) {\tiny $61(v^5)$};
         \node[above right] at (1,4) {\tiny $63(v^5)$};
         \node[label={[above right, yshift=1ex]:\tiny $\un{66(v^6)}$}] at (1,4) {};      
         
         \node[above left] at (1,10) {\tiny $\un{54(v^4)}$};
         \node[label={[above left, yshift=1ex]:\tiny $57(v^5)$}] at (1,10) {};
         \node[above left] at (2,10) {\tiny $59(v^5)$};
         \node[above left] at (3,10) {\tiny $61(v^5)$};
         \node[above left] at (4,10) {\tiny $63(v^5)$};
         \node[label={[above left, yshift=1ex]:\tiny $\un{66(v^6)}$}] at (4,10) {};
         
         \node[below right] at (3,5) {\tiny $\un{77(v^7)}$};
         \node[below left] at (2,13) {\tiny $\un{77(v^7)}$};
         \node[right] at (5,7) {\tiny $\un{77(v^7)}$};
         
         \node[above right] at (3,7) {\color{blue}\tiny $\un{88(v^8)}$};
         \node[above left] at (2,15) {\tiny $\un{88(v^8)}$};
         
         \node[below right] at (1,10) {\tiny $\un{99(v^9)}$};
         \node[label={[below right, yshift=-3ex]:\tiny $102(v^{10})$}] at (1,10) {};
         \node[below right] at (2,10) {\tiny $104(v^{10})$};
         \node[below right] at (3,10) {\tiny $106(v^{10})$};
         \node[below right, xshift=-1ex] at (4,10) {\tiny $108(v^{10})$};
         \node[label={[below right, xshift=-1ex, yshift=-2.5ex]:\tiny $\un{111(v^{11})}$}] at (4,10) {};
         
         \node[above right] at (5,14) {\tiny $\un{99(v^9)}$};
         \node[right] at (5,14) {\tiny $102(v^{10})$};
         \node[right] at (5,13) {\tiny $104(v^{10})$};
         \node[right] at (5,12) {\tiny $106(v^{10})$};
         \node[above right] at (5,11) {\tiny $108(v^{10})$};
         \node[right] at (5,11) {\tiny $\un{111(v^{11})}$};
         
         \node[below left] at (4,16) {\tiny $\un{99(v^9)}$};
         \node[label={[below left, yshift=-3ex]:\tiny $102(v^{10})$}] at (4,16) {};
         \node[below left] at (3,17) {\tiny $104(v^{10})$};
         \node[below left] at (2,18) {\tiny $106(v^{10})$};
         \node[below left, xshift=1ex] at (1,19) {\tiny $108(v^{10})$};
         \node[label={[below left, xshift=1ex, yshift=-2.5ex]:\tiny $\un{111(v^{11})}$}] at (1,19) {};
         
         \node[above right] at (2,13) {\color{blue}\tiny $\un{122(v^{12})}$};
         \node[above left] at (2,20) {\tiny $\un{122(v^{12})}$};
      \end{tikzpicture}}

      \caption{$(5\varpi_1, 10(v^0))$}
      \label{figl5k1}
   \end{minipage}
   \hfill
   \begin{minipage}[b]{.64\textwidth}
      \centering
      \resizebox{\textwidth}{!}{
      \begin{tikzpicture}[scale=1.1, rotate=30, point/.style={inner sep=0pt,outer sep = 1mm,circle,fill,minimum size=1mm}]
         \pgftransformcm{1}{0}{0.5}{sqrt(3)/2}{\pgfpoint{0cm}{0cm}};   
         \draw[opacity=0.5, ->] (0,0) -- (11,0);
         \draw[opacity=0.5, ->] (0,0) -- (0,22);
         \foreach \y in {0,1, ..., 20} 
         {
            \node (x5-y\y) at (5,\y) [point] {};
            \ifthenelse{\y > 0}{
               \pgfmathtruncatemacro{\yprev}{\y - 1}
               \draw[color=gray!70,->] (x5-y\y) -- (x5-y\yprev);}
               {}
         }
         \foreach \y in {0,1, ..., 20}
         {
            \node (x10-y\y) at (10,\y) [point] {};
            \ifthenelse{\y > 0}{
               \pgfmathtruncatemacro{\yprev}{\y - 1}
               \draw[color=gray!70,->] (x10-y\y) -- (x10-y\yprev);}
               {}
         }      
         \foreach \t in {1, ..., 5}
         {
            \pgfmathtruncatemacro{\x}{5 - \t}
            \node (x\x-y\t) at (\x,\t) [point] {};         
            \pgfmathtruncatemacro{\xprev}{\x+1}
            \pgfmathtruncatemacro{\yprev}{\t-1}
            \draw[color=gray!70,->] (x\xprev-y\yprev) -- (x\x-y\t);
            
            \pgfmathtruncatemacro{\x}{10 - \t}
            \pgfmathtruncatemacro{\y}{15 + \t}
            \node (x\x-y\y) at (\x,\y) [point] {};         
            \pgfmathtruncatemacro{\xprev}{\x+1}
            \pgfmathtruncatemacro{\yprev}{\y-1}
            \draw[color=gray!70,->] (x\xprev-y\yprev) -- (x\x-y\y);
         }
         \foreach \ty in {0, 5, 10} {    
            \foreach \tx in {1, ..., 10}
            {
               \pgfmathtruncatemacro{\x}{10 - \tx}
               \pgfmathtruncatemacro{\y}{\ty + \tx}
               \node (x\x-y\y) at (\x,\y) [point] {};
               \pgfmathtruncatemacro{\xprev}{\x+1}
               \pgfmathtruncatemacro{\yprev}{\y-1}
               \draw[color=gray!70,->] (x\xprev-y\yprev) -- (x\x-y\y);
            }
         }
         \foreach \y in {5, 10, 15, 20} {
            \foreach \x in {1, ..., 10} 
            {
               \node (x\x-y\y) at (\x,\y) [point] {};
               \pgfmathtruncatemacro{\xprev}{\x-1}
               \draw[color=gray!70,->] (x\xprev-y\y) -- (x\x-y\y);
            }
         }
         \node[below right] at (10,0) {\tiny $\un{20(v^0)}$};
         \node[below left] at (9,1) {\tiny $22(v^0)$};
         \node[below left] at (8,2) {\tiny $24(v^0)$};
         \node[below left] at (7,3) {\tiny $26(v^0)$};
         \node[below left] at (6,4) {\tiny $28(v^0)$};
         \node[left] at (6,4) {\tiny $\un{31(v^1)}$};
         \node[right] at (5,2) {\tiny $\un{42(v^2)}$};
         \node[above left] at (8,5) {\tiny $\un{42(v^2)}$};      
         \node[below left] at (3,2) {\tiny $\un{53(v^3)}$};
         \node[right] at (10,3) {\tiny $\un{53(v^3)}$};
         \node[below left] at (8,7) {\tiny $\un{53(v^3)}$};
         
         \node[above left] at (1,5) {\tiny $\un{64(v^4)}$};
         \node[label={[above left, yshift=1ex]:\tiny $67(v^5)$}] at (1,5) {};
         \node[above left] at (2,5) {\tiny $69(v^5)$};
         \node[above left] at (3,5) {\tiny $71(v^5)$};
         \node[above left] at (4,5) {\tiny $73(v^5)$};
         \node[label={[above left, yshift=1ex]:\tiny $\un{76(v^6)}$}] at (4,5) {};
         
         \node[above right] at (9,1) {\tiny $\un{64(v^4)}$};
         \node[label={[above right, yshift=1ex]:\tiny $67(v^5)$}] at (9,1) {};
         \node[above right] at (8,2) {\tiny $69(v^5)$};
         \node[above right] at (7,3) {\tiny $71(v^5)$};
         \node[above right] at (6,4) {\tiny $73(v^5)$};
         \node[label={[above right, yshift=1ex]:\tiny $\un{76(v^6)}$}] at (6,4) {};      
         
         \node[above left] at (6,10) {\tiny $\un{64(v^4)}$};
         \node[label={[above left, yshift=1ex]:\tiny $67(v^5)$}] at (6,10) {};
         \node[above left] at (7,10) {\tiny $69(v^5)$};
         \node[above left] at (8,10) {\tiny $71(v^5)$};
         \node[above left] at (9,10) {\tiny $73(v^5)$};
         \node[label={[above left, yshift=1ex]:\tiny $\un{76(v^6)}$}] at (9,10) {};
         
         \node[above left] at (5,9) {\tiny $\un{64(v^4)}$};
         \node[left] at (5,9) {\tiny $67(v^5)$};
         \node[left] at (5,8) {\tiny $69(v^5)$};
         \node[left] at (5,7) {\tiny $71(v^5)$};
         \node[above left] at (5,6) {\tiny $73(v^5)$};
         \node[left] at (5,6) {\tiny $\un{76(v^6)}$};     
         
         \node[left] at (5,2) {\color{red}\tiny $\un{87(v^7)}$};
         \node[below left] at (2,8) {\color{red}\tiny $\un{87(v^7)}$};
         \node[below right] at (8,5) {\color{red}\tiny $\un{87(v^7)}$};
         \node[below left] at (7,13) {\tiny $\un{87(v^7)}$};
         \node[right] at (10,7) {\tiny $\un{87(v^7)}$};
         
         \node[above right] at (3,2) {\tiny $\un{98(v^8)}$};
         \node[above left] at (2,10) {\tiny $\un{98(v^8)}$};
         \node[above right] at (8,7) {\color{blue}\tiny $\un{98(v^8)}$};
         \node[above left] at (7,15) {\tiny $\un{98(v^8)}$};
         \node[left] at (5,13) {\tiny $\un{98(v^8)}$};
         
         \node[below right] at (1,5) {\tiny $\un{109(v^9)}$};
         \node[label={[below right, yshift=-3ex]:\tiny $112(v^{10})$}] at (1,5) {};
         \node[below right] at (2,5) {\tiny $114(v^{10})$};
         \node[below right] at (3,5) {\tiny $116(v^{10})$};
         \node[below right, xshift=-1ex] at (4,5) {\tiny $118(v^{10})$};
         \node[label={[below right, xshift=-1ex, yshift=-2.5ex]:\tiny $\un{121(v^{11})}$}] at (4,5) {};
         
         \node[below right] at (6,10) {\color{red}\tiny $\un{109(v^9)}$};
         \node[label={[below right, yshift=-3ex]:\tiny $112(v^{10})$}] at (6,10) {};
         \node[below right] at (7,10) {\tiny $114(v^{10})$};
         \node[below right] at (8,10) {\tiny $116(v^{10})$};
         \node[below right, xshift=-1ex] at (9,10) {\tiny $118(v^{10})$};
         \node[label={[below right, xshift=-1ex, yshift=-2.5ex]:\tiny $\un{121(v^{11})}$}] at (9,10) {};
         
         \node[above right] at (5,9) {\color{red}\tiny $\un{109(v^9)}$};
         \node[right] at (5,9) {\tiny $112(v^{10})$};
         \node[right] at (5,8) {\tiny $114(v^{10})$};
         \node[right] at (5,7) {\tiny $116(v^{10})$};
         \node[above right] at (5,6) {\tiny $118(v^{10})$};
         \node[right] at (5,6) {\tiny $\un{121(v^{11})}$};
         
         \node[below left] at (4,11) {\color{red}\tiny $\un{109(v^9)}$};
         \node[label={[below left, yshift=-3ex]:\tiny $112(v^{10})$}] at (4,11) {};
         \node[below left] at (3,12) {\tiny $114(v^{10})$};
         \node[below left] at (2,13) {\tiny $116(v^{10})$};
         \node[below left, xshift=1ex] at (1,14) {\tiny $118(v^{10})$};
         \node[label={[below left, xshift=1ex, yshift=-2.5ex]:\tiny $\un{121(v^{11})}$}] at (1,14) {};
         
         \node[above right] at (10,14) {\tiny $\un{109(v^9)}$};
         \node[right] at (10,14) {\tiny $112(v^{10})$};
         \node[right] at (10,13) {\tiny $114(v^{10})$};
         \node[right] at (10,12) {\tiny $116(v^{10})$};
         \node[above right] at (10,11) {\tiny $118(v^{10})$};
         \node[right] at (10,11) {\tiny $\un{121(v^{11})}$};
         
         \node[below left] at (9,16) {\tiny $\un{109(v^9)}$};
         \node[label={[below left, yshift=-3ex]:\tiny $112(v^{10})$}] at (9,16) {};
         \node[below left] at (8,17) {\tiny $114(v^{10})$};
         \node[below left] at (7,18) {\tiny $116(v^{10})$};
         \node[below left, xshift=1ex] at (6,19) {\tiny $118(v^{10})$};
         \node[label={[below left, xshift=1ex, yshift=-2.5ex]:\tiny $\un{121(v^{11})}$}] at (6,19) {};
         
         \node[above right] at (2,8) {\color{blue}\tiny $\un{132(v^{12})}$};
         \node[above right] at (7,13) {\color{blue}\tiny $\un{132(v^{12})}$};
         \node[above left] at (3,15) {\tiny $\un{132(v^{12})}$};
         \node[right] at (5,17) {\tiny $\un{132(v^{12})}$};
         \node[above left] at (7,20) {\tiny $\un{132(v^{12})}$};
      \end{tikzpicture}}

      \caption{$(10\varpi_1, 20(v^0))$}
      \label{figl5k2}
   \end{minipage}\\
   \vfill
   These figures show the multiset for $\ell = 5$ after the first $10$ iterations 
   applied to
\end{figure*}

\section{The current Data}
\label{secData}

\begin{sidewaysfigure}
   \includegraphics[width=\textwidth]{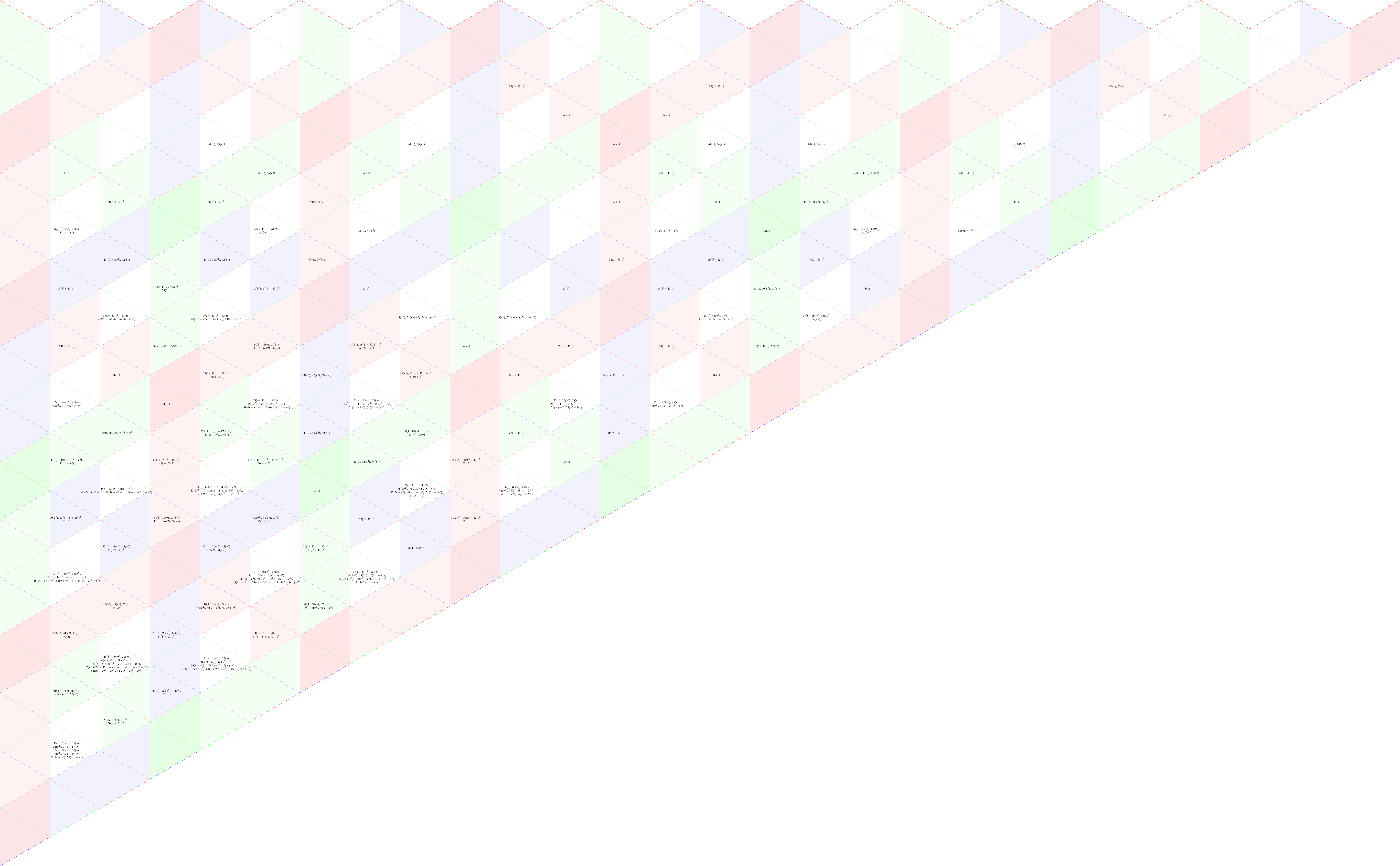}
   \captionof{figure}{Triples for $p=3$}
   \label{figp3}
\end{sidewaysfigure}

\begin{sidewaysfigure}
   \includegraphics[width=\textwidth]{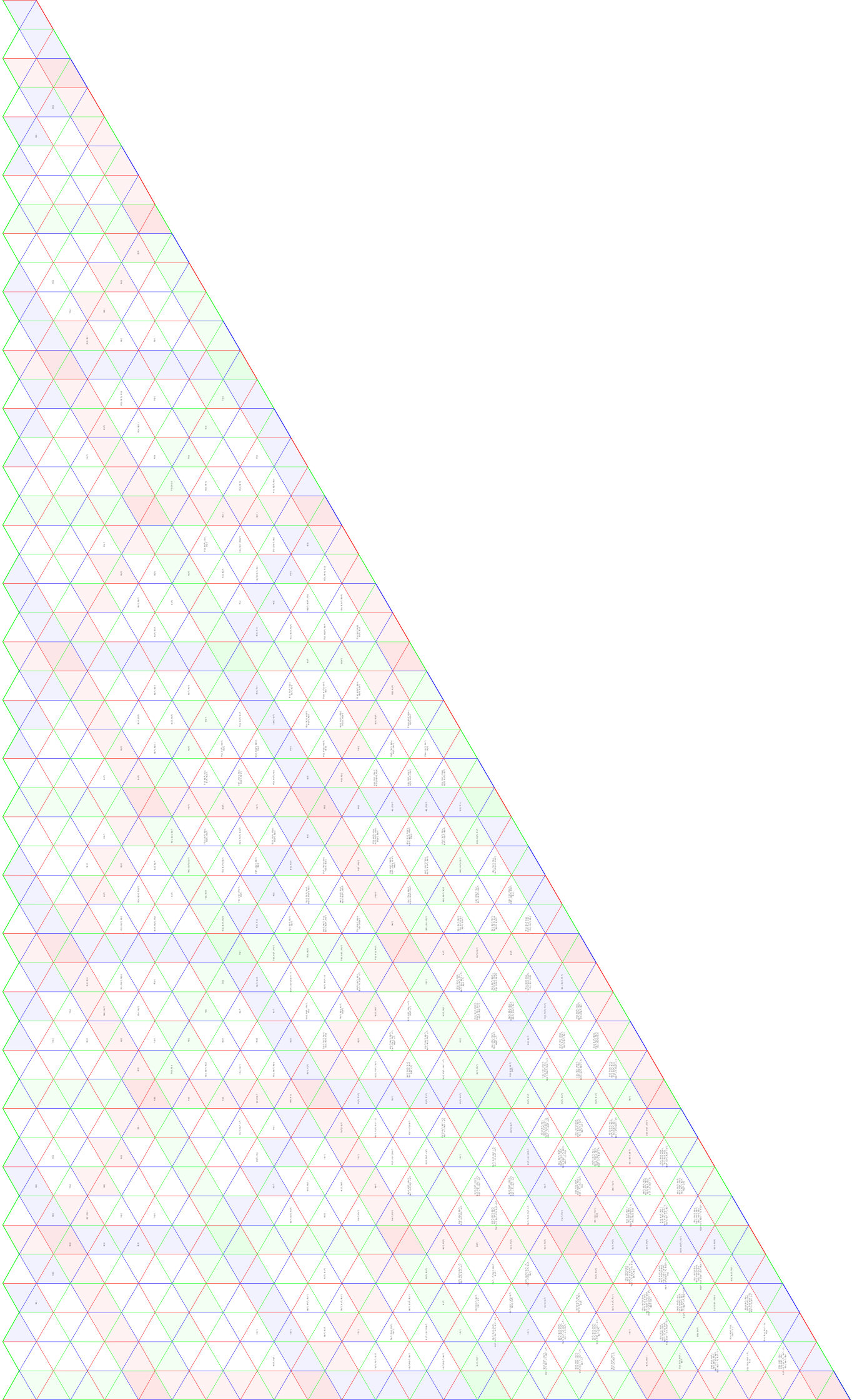}
   \captionof{figure}{Triples for $p=5$}
   \label{figp5}
\end{sidewaysfigure}

\begin{sidewaysfigure}
   \includegraphics[width=\textwidth]{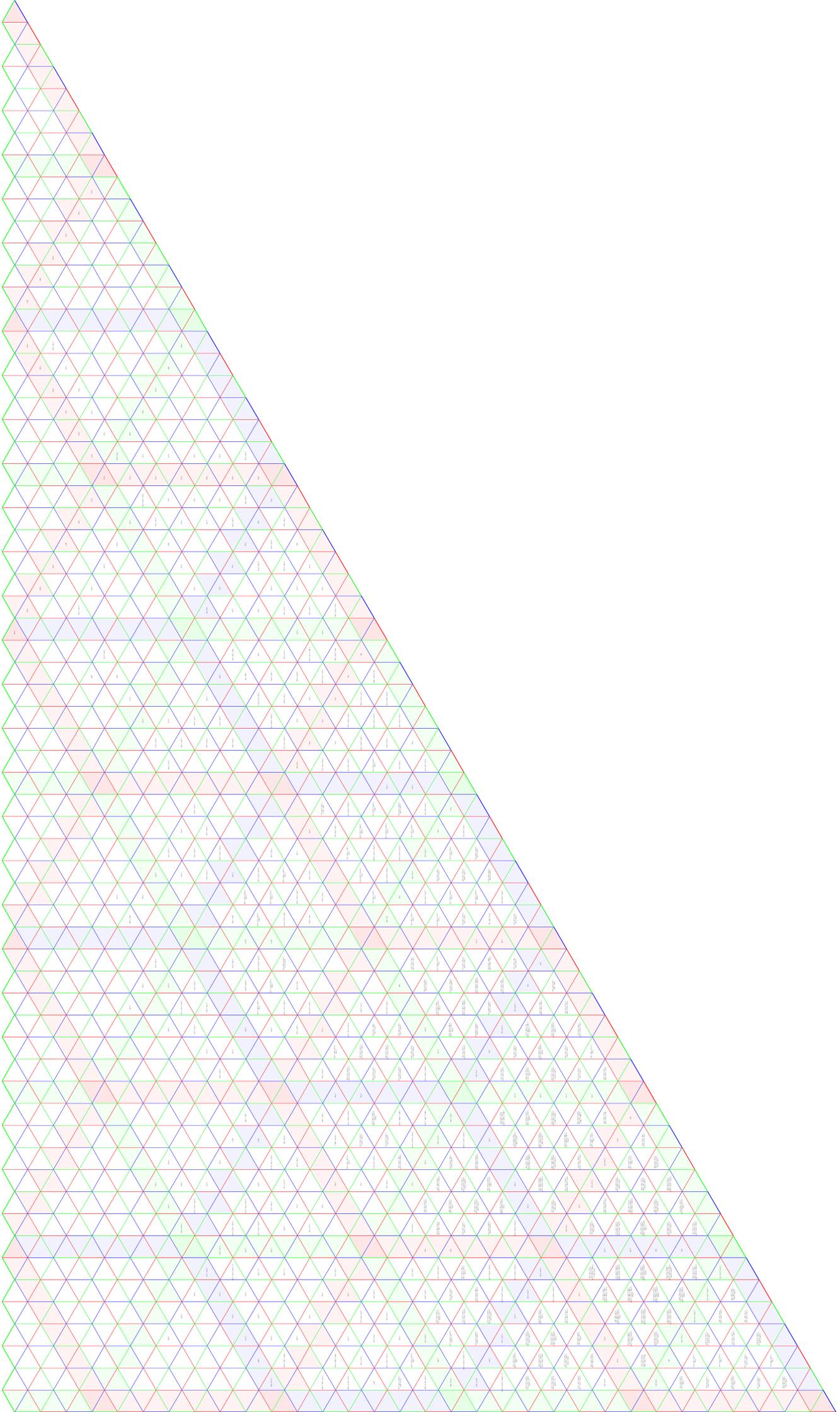}
   \captionof{figure}{Triples for $p=7$}
   \label{figp7}
\end{sidewaysfigure}

\begin{sidewaysfigure}
   \includegraphics[width=\textwidth]{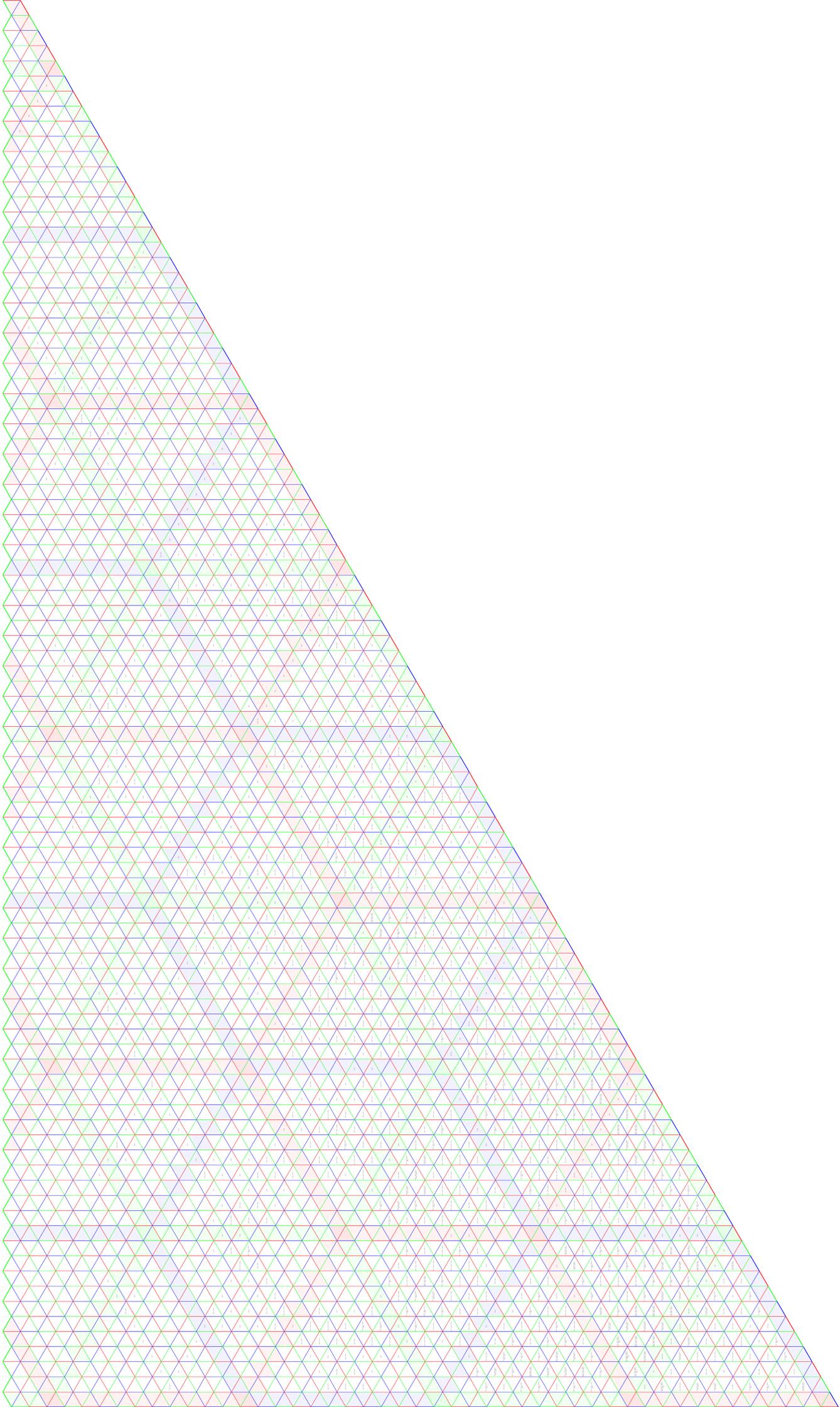}
   \captionof{figure}{Triples for $p=11$}
   \label{figp11}
\end{sidewaysfigure}

The current data for \hyperref[figp3]{$\ell=3$}, \hyperref[figp5]{$\ell=5$}, 
\hyperref[figp7]{$\ell=7$} and \hyperref[figp11]{$\ell=11$} is displayed in four 
figures at the end of the paper.

\subsection{Preparation of the Data}
\label{secPrep}

Imagine we have calculated the elements $\pkl{x_i}$ for $1 \leqslant i \leqslant n$.
In this section we will explain how we have obtained the figures at the end of 
the paper from the given data. We have applied the following steps to process the data:

\begin{enumerate}
   \item Express the elements $\pkl{x_i}$ for $1 \leqslant i \leqslant n$
         in terms of the Kazhdan-Lusztig basis $\{ \kl{x} \; \vert \; x \in \fW \}$
         of the anti-spherical module.
   \item Illustrate all calculated $p$-Kazhdan-Lusztig basis elements
         in one combined picture as follows:
         For a Laurent polynomial $f$ let $\bar{f}$ be the polynomial 
         obtained by forgetting all negative powers of $v$ and their 
         coefficients in $f$.
         If $\kl{y}$ occurs with coefficient $f$ in $\pkl{x_j}$, then
         write $j(\bar{f})$ in the $p$-alcove corresponding to $y$.
   \item For all $\mu \in \charlat_{++}$ replace each ``triple'' of the form
         \[
         \begin{tikzpicture}[point/.style={inner sep=0pt,outer sep = 1mm,circle,fill,minimum size=1mm}]
            \pgftransformcm{1}{0}{0.5}{sqrt(3)/2}{\pgfpoint{0cm}{0cm}};
            \draw (0,0) -- (0,3) -- (3,3) -- (3,0) -- cycle;
            \draw (0,3) -- (3,0);
            \node[labelling] at (2/3,5/3) {\footnotesize$i(f),$};
            \node[labelling] at (1, 1) {\footnotesize$(i+2)(f)$};
            \node[labelling] at (2, 2) {\footnotesize$(i+1)(f)$};
         \end{tikzpicture} \]
         by
         \[
         \begin{tikzpicture}[point/.style={inner sep=0pt,outer sep = 1mm,circle,fill,minimum size=1mm}]
            \pgftransformcm{1}{0}{0.5}{sqrt(3)/2}{\pgfpoint{0cm}{0cm}};
            \draw (0,0) -- (0,3) -- (3,3) -- (3,0) -- cycle;
            \draw (0,3) -- (3,0);
            \node[labelling] at (1, 1) {\footnotesize$i(f)$};
         \end{tikzpicture} \text{.}\]
         These triples are consequences of the Kazhdan-Lusztig star operations
         as follows: 
         Denote the simple reflections as $s \defeq x_i^{-1} x_{i+1}$, 
         $t \defeq x_{i+1}^{-1} x_{i+2}$ and $u \defeq x_{i+2}^{-1} x_{i+3}$.
         Observe that the wall separating the two $p$-alcoves in the box
         $B_{\mu}$ will be $s$-coloured. So one may apply 
         \cite[Proposition 4.5]{JeABC} first to the two rank $2$ standard 
         parabolic subgroup given by $\{s, u\}$ and then to $\{s, t\}$ to
         obtain the triple shown above.\footnote{For illustration see
         \url{https://www.maths.usyd.edu.au/u/geordie/pCanA2/p5pretriples.pdf}.}
\end{enumerate}
The resulting figures are depicted at the end of the paper.

It should be noted that in order to compare the figures with the predictions
of the conjecture, the reader still needs to use the following heuristic to
remove third generation phenomena. Even though there is no precise definition
of $\{ \pklgen{x}{3} \; \vert \; x \in \fW \}$ in general, we have a clear
understanding of what it is for $SL_2$ (see explanation and \cref{figSL2pCan}
below). 
\begin{heuristic}
   Remove any labelled point $(\mu, n(v^k))$ (and all labelled
   points induced by it) whose restriction to a rank $1$ Levi subgroup 
   gives a third generation contribution.
\end{heuristic}

For $SL_2$ the affine Weyl group is isomorphic to the infinite dihedral group
\[ W = \langle s, t \; \vert \; s^2 = t^2 = \id \rangle\text{.} \]
Let ${}_s \hat{k}$ denote the alternating word in $s$ and $t$ of length $k$
starting in $s$. Suppose that $t$ is the simple reflection generating the 
finite Weyl group corresponding to our chosen maximal torus. For $p=3$ we get for
the $p$-Kazhdan-Lusztig basis of the anti-spherical module:
\begin{alignat*}{9}
   &\pkl[3]{s}        &&= \kl{s} && && && && && && &&\\
   &\pkl[3]{st}       &&=        &&\kl{st} && && && && && &&\\
   &\pkl[3]{sts}      &&=        &&        &&\kl{sts} && && && && &&\\
   &\pkl[3]{stst}     &&=        &&\kl{st} &&\ \  +   &&\kl{stst} && && && &&\\
   &\pkl[3]{ststs}    &&= \kl{s} &&        &&\ \  +   &&          &&\kl{ststs} && && &&\\
   &\pkl[3]{ststst}   &&=        &&        &&         &&          &&           &&\kl{ststst} && &&\\
   &\pkl[3]{stststs}  &&=        &&        &&         &&          &&\kl{ststs} &&\ \ \ +       &&\kl{stststs} &&\\
   &\pkl[3]{stststst} &&=        &&        &&         &&\kl{stst} &&           &&\ \ \ +       &&             &&\kl{stststst} \\
\end{alignat*}
For more information on how to calculate the $p$-Kazhdan-Lusztig basis
and its properties in this case, we refer the reader to \cite[\S9.1 and \S7.2]{JpKLTheory}.
For $SL_2$, the base change coefficients between the $p$-Kahzdan-Lusztig basis
and the Kazhdan-Lusztig basis all are either $0$ or $1$. In \cref{figSL2pCan} we 
have encoded whenever $\kl{ {}_s\hat{m}}$ occurs with non-trivial coefficient
in $\pkl{{}_s \hat{n}}$ for $p=3$ by a coloured box. Higher generation contributions 
are coloured in increasingly lighter shades of grey. Overall the figure shows
contributions from the first four generations.

\def \xspace{0.5}
\def \yspace{0.5}
\def \num{41}

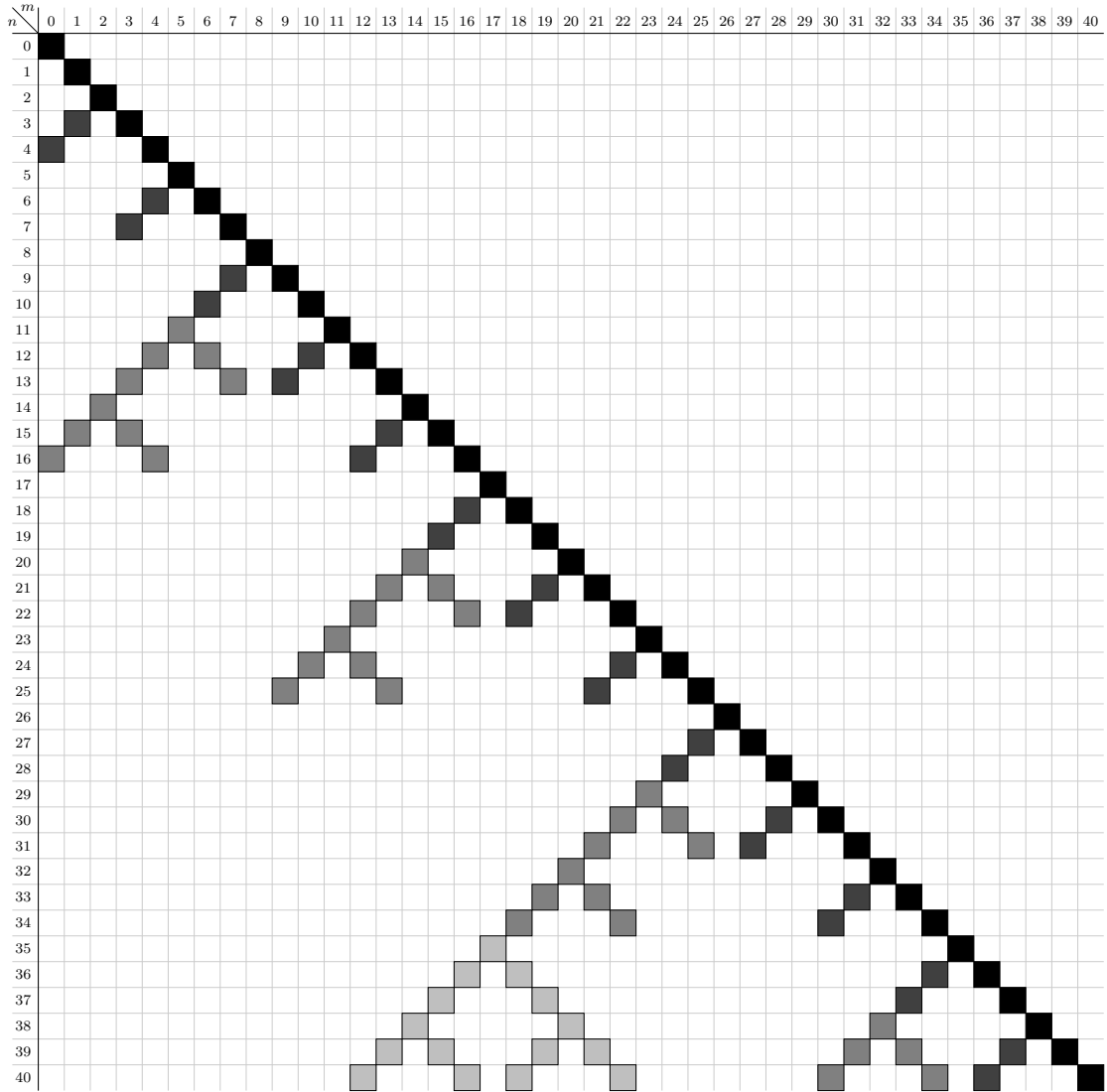
\begin{figure}[htb]
   \centering
   \resizebox{\textwidth}{!}{
   \begin{tikzpicture}[auto, baseline=(current  bounding  box.center)]
      \pgfmathsetmacro{\max}{\num-1};
      \foreach \i in {0, 1, ..., \max} {
         \pgfmathsetmacro{\x}{\i*\xspace + \xspace/2};
         \node[labelling, label=above:\footnotesize$\i$] (\i) at (\x,0) {};
         \node[labelling, label=left:\footnotesize$\i$] (T\i) at (0,-\x) {};
         \draw[color=black!20] (\i*\xspace,\yspace) to (\i*\xspace, -\num*\yspace);
         \draw[color=black!20] (-\xspace, -\i*\yspace) to (\num*\xspace,-\i*\yspace);
      }
      
      \foreach \gen/\x/\y in {1/1/3, 1/0/4, 1/4/6, 1/3/7, 1/7/9, 1/6/10, 2/5/11, 2/4/12, 2/6/12, 2/3/13, 2/7/13, 
                         2/2/14, 2/1/15, 2/3/15, 2/0/16, 2/4/16, 1/10/12, 1/9/13, 1/13/15, 1/12/16,
                         1/16/18, 1/15/19, 2/14/20, 2/13/21, 2/15/21, 2/12/22, 2/16/22, 2/11/23, 2/10/24, 2/12/24, 2/9/25, 2/13/25,
                         1/19/21, 1/18/22, 1/22/24, 1/21/25, 1/25/27, 1/24/28, 2/23/29, 2/22/30, 2/24/30, 2/21/31, 2/25/31,
                         2/20/32, 2/19/33, 2/21/33, 2/18/34, 2/22/34, 3/17/35, 3/16/36, 3/18/36, 3/15/37, 3/19/37, 3/14/38, 3/20/38, 3/13/39, 3/15/39, 3/19/39, 3/21/39,
                         1/28/30, 1/27/31, 1/31/33, 1/30/34, 1/34/36, 1/33/37, 2/32/38, 2/31/39, 2/33/39, 1/37/39,
                         3/12/40, 3/16/40, 3/18/40, 3/22/40, 2/30/40, 2/34/40, 1/36/40} {
         \pgfmathsetmacro{\Px}{\x*\xspace}
         \pgfmathsetmacro{\Py}{-\y*\yspace}
         \pgfmathsetmacro{\Qx}{(\x + 1)*\xspace}
         \pgfmathsetmacro{\Qy}{-(\y + 1)*\yspace}
         \pgfmathsetmacro{\shade}{100 - (\gen * 25)}
         \draw[fill=black!\shade] (\Px, \Py) rectangle (\Qx, \Qy);
      }
      
      \foreach \i in {0, 1, ..., \max} {
         \pgfmathsetmacro{\Px}{\i*\xspace}
         \pgfmathsetmacro{\Py}{-\i*\yspace}
         \pgfmathsetmacro{\Qx}{(\i + 1)*\xspace}
         \pgfmathsetmacro{\Qy}{-(\i + 1)*\yspace}
         \draw[fill=black] (\Px, \Py) rectangle (\Qx, \Qy);
      }
      \draw (-\xspace,0) to (\num*\xspace,0);
      \draw (0,\yspace) to (0, -\num*\yspace);
      \draw (0,0) to (-\xspace, \yspace);
      \node[labelling, label={center:\footnotesize$n$}] (n) at (-\xspace, 0.2cm) {};
      \node[labelling, label={center:\footnotesize$m$}] (m) at (-0.2cm, \yspace) {};
   \end{tikzpicture}}
   \caption{Generations for the $p$-Kazhdan-Lusztig basis for $SL_2$ for $p=3$.}
   \label{figSL2pCan}
\end{figure}

\subsection{Examples for the Merging Rule}

In total, there are ten examples for the merging rule mentioned above in 
the current data. The reader is invited to locate each one of them in the data.
\begin{description}
   \item[$\ell = 5$]: Merging of type \RN{2} produces $88(v^8)$ around $\ell\varpi_1 + \ell\varpi_2$.\\
                   Merging of type \RN{3} produces $87(v^7)$ around $\ell\varpi_1 + \ell\varpi_2$.
   \item[$\ell = 7$]: Merging of type \RN{2} produces $104(v^6)$ around $\ell\varpi_1 + \ell\varpi_2$.\\
                   Merging of type \RN{3} produces $103(v^5)$ around $\ell\varpi_1 + \ell\varpi_2$.\\
                   Merging of type \RN{2} produces $118(v^6)$ around $2\ell\varpi_1 + \ell\varpi_2$.\\
                   Merging of type \RN{3} produces $117(v^5)$ around $2\ell\varpi_1 + \ell\varpi_2$.
   \item[$\ell = 11$]: Merging of type \RN{2} produces $160(v^6)$ around $\ell\varpi_1 + \ell\varpi_2$.\\
                    Merging of type \RN{3} produces $159(v^5)$ around $\ell\varpi_1 + \ell\varpi_2$.\\
                    Merging of type \RN{2} produces $182(v^6)$ around $2\ell\varpi_1 + \ell\varpi_2$.\\
                    Merging of type \RN{3} produces $181(v^5)$ around $2\ell\varpi_1 + \ell\varpi_2$.
\end{description}

\printbibliography

\Address

\newgeometry{scale=1}
\thispagestyle{empty}



\restoregeometry

\end{document}